\documentclass[a4paper,12pt]{amsart}
\usepackage{amsfonts,amsmath,amsthm,amssymb,multicol}
\usepackage{enumitem}
\usepackage{mathrsfs}
\usepackage[margin=1.15in]{geometry}
\usepackage[colorlinks=true,linkcolor=blue,filecolor=blue,citecolor=red]{hyperref}
\usepackage{orcidlink}
\usepackage{cleveref}

\theoremstyle{plain}
\newtheorem{theorem}{Theorem}[section]

\newtheorem{conjecture}[theorem]{Conjecture}
\theoremstyle{definition}

\newtheorem{remark}[theorem]{Remark}

\begin{document}
\title[The $k$-elongated plane partition function modulo small powers of $5$]{The $k$-elongated plane partition function modulo small powers of $5$}
\author[Russelle Guadalupe]{Russelle Guadalupe\orcidlink{0009-0001-8974-4502}}
\address{Institute of Mathematics, University of the Philippines, Diliman\\
Quezon City 1101, Philippines}
\email{rguadalupe@math.upd.edu.ph}

\renewcommand{\thefootnote}{}

\footnote{2020 \emph{Mathematics Subject Classification}: Primary 11P83, Secondary 05A17}

\footnote{\emph{Key words and phrases}: plane partitions, Ramanujan-type congruences, $q$-series, dissection formulas}

\renewcommand{\thefootnote}{\arabic{footnote}}

\setcounter{footnote}{0}

\begin{abstract}
Andrews and Paule revisited combinatorial structures known as the $k$-elongated partition diamonds, which were introduced in connection with the study of the broken $k$-diamond partitions. They found the generating function for the number $d_k(n)$ of partitions obtained by summing the links of such partition diamonds of length $n$ and discovered congruences for $d_k(n)$ using modular forms. Since then, congruences for $d_k(n)$ modulo certain powers of primes have been proven via elementary means and modular forms by many authors, most recently Banerjee and Smoot who established an infinite family of congruences for $d_5(n)$ modulo powers of $5$. We extend in this paper the list of known results for $d_k(n)$ by proving infinite families of congruences for $d_k(n)$ modulo $5,25$, and $125$ using classical $q$-series manipulations and $5$-dissections.
\end{abstract}

\maketitle

\section{Introduction}\label{sec1}

Throughout this paper, define $f_m := \prod_{n\geq 1}(1-q^{mn})$ for a positive integer $m$ and a complex number $q$ with $|q| < 1$. In 2022, Andrews and Paule \cite{andpau2} revisited certain combinatorial objects called the $k$-elongated partition diamonds, which were defined fifteen years ago in relation with the study of the broken $k$-diamond partitions \cite{andpau1}. They considered these partition diamonds as one of the examples of Schmidt type partitions arising from partitions on certain graphs. Using MacMahon's partition analysis, Andrews and Paule \cite{andpau2} obtained the following generating function 
\begin{align*}
\sum_{n\geq 0} d_k(n)q^n = \dfrac{f_2^k}{f_1^{3k+1}}
\end{align*}
for the number $d_k(n)$ of partitions found by summing the links of the $k$-elongated partition diamonds of length $n$. They then discovered several congruences for $d_1, d_2$, and $d_3$ modulo certain powers of primes by primarily using the Mathematica package \texttt{RaduRK} implemented by Smoot \cite{smoot1}, which is based on the Ramanujan-Kolberg algorithm presented by Radu \cite{radu2}.

Since the inception of the function $d_k(n)$, various authors have extended its congruence properties through elementary $q$-series manipulations and modular forms. da Silva, Hirschhorn, and Sellers \cite{dshs} gave elementary proofs of several congruences for $d_k(n)$ found by Andrews and Paule \cite{andpau2} and added new individual congruences modulo small prime powers. Yao \cite{yao} provided elementary proofs of the congruences modulo $81, 243$, and $729$ for $d_k(n)$ conjectured by Andrews and Paule \cite[Conjectures 1 and 2]{andpau2}. Baruah, Das, and Talukdar \cite{bdt} found infinite families of congruences for $d_k(n)$ modulo powers of $2$ and $3$ and proved the following refinement of the result of da Silva, Hirschhorn, and Sellers \cite[Theorem 4.1]{dshs} on the existence of infinite congruence families for $d_k(n)$ modulo prime powers. 

\begin{theorem}{\cite[Theorem 5.1]{bdt}}\label{thm11}
Let $p$ be a prime, $k\geq 1, N\geq 1, M\geq 1$, and $r$ be integers such that $1\leq r\leq p^M-1$. If
\begin{align*}
d_k(p^Mn+r)\equiv 0\pmod{p^N}
\end{align*}
for all $n\geq 0$, then
\begin{align*}
d_{p^{M+N-1}j+k}(p^Mn+r)\equiv 0\pmod{p^N}
\end{align*}
for all $n\geq 0$ and $j\geq 0$.
\end{theorem}

Recently, Banerjee and Smoot \cite{bansmo} applied the localization method and modular cusp analysis to establish the following congruence 
\begin{align}\label{eq11}
d_5(n)\equiv 0\pmod{5^{\lfloor k/2\rfloor+1}}
\end{align}
for all $k\geq 1$ and $n\geq 1$ such that $4n\equiv 1\pmod{5^k}$. Prior to this remarkable result, Smoot \cite{smoot2} also applied the same method to derive the following refinement of the 
conjectural congruence of Andrews and Paule \cite[Conjecture 3]{andpau2} given by
\begin{align*}
d_2(n)\equiv 0\pmod{3^{2\lfloor k/2\rfloor+1}}
\end{align*}
for all $k\geq 1$ and $n\geq 1$ such that $8n\equiv 1\pmod{3^k}$. We note that these congruences satisfied by $d_2$ and $d_5$ resemble the groundbreaking result of Ramanujan \cite{berno}, Watson \cite{wat}, and Atkin \cite{atk} on the congruences modulo powers of $\ell\in \{5,7,11\}$ for the number $p(n)$ of unrestricted partitions of $n$ given by
\begin{align}\label{eq12}
p(n)\equiv 0\pmod{\ell^\beta}
\end{align}
for all $\alpha\geq 1$ and $n\geq 1$ such that $24n\equiv 1\pmod{\ell^\alpha}$, where $\beta := \alpha$ if $\ell\in \{5,11\}$ and $\beta := \lfloor\alpha/2\rfloor+1$ if $\ell=7$.

The goal of this paper is to extend the list of the known congruences for $d_k(n)$ by deriving infinite families of congruences modulo $5, 25$, and $125$. More precisely, we aim to prove our main results as shown by purely elementary methods. 

\begin{theorem}\label{thm12}
For all $c\geq 0$ and $n\geq 0$, 
\begin{align*}
d_{25c+16}(25n+8)&\equiv 0\pmod{5},\\
d_{25c+8}(5n+1)&\equiv 0\pmod{25},\\
d_{25c+k}(25n+24-k)&\equiv 0\pmod{25}
\end{align*}
for $k\in\{0,5,10\}$.
\end{theorem}

\begin{remark}\label{rem13}
\begin{enumerate}
\item The first congruence of Theorem \ref{thm12} slightly improves the recent congruence of Banerjee and Smoot \cite[Theorem 1.3]{bansmo} given by 
\begin{align*}
d_{75c+16}(25n+8)\equiv 0\pmod{5}
\end{align*}
for all $c\geq 0$ and $n\geq 0$.
\item Using Theorem \ref{thm11} and \cite[Theorem 1.5]{bansmo}, we infer that
\begin{align*}
d_{25c+k}(25n+24-k)\equiv 0\pmod{5}
\end{align*}
for $k\in \{1,3,13,15,18,20,23\}$ and for all $c\geq 0$ and $n\geq 0$.
\end{enumerate}
\end{remark}

\begin{theorem}\label{thm14}	
For all $c\geq 0$ and $n\geq 0$, 
\begin{align*}
&
\begin{aligned}[t]
d_{125c+6}(125n+43,93)&\equiv 0\pmod{5},\\
d_{125c+61}(125n+38,63,88,113)&\equiv 0\pmod{5},\\
d_{125c+82}(125n+67,92,117)&\equiv 0\pmod{5},\\
d_{125c+106}(125n+68,118)&\equiv 0\pmod{5},\\
d_{125c+111}(125n+63,113)&\equiv 0\pmod{5},\\
d_{125c+46}(125n+28,78,103)&\equiv 0\pmod{5},\\
d_{125c+71}(125n+53,78)&\equiv 0\pmod{5},\\
d_{125c+96}(125n+103)&\equiv 0\pmod{5},\\
d_{125c+2}(125n+97,122)&\equiv 0\pmod{5},\\
d_{125c+107}(125n+42,92)&\equiv 0\pmod{5},\\
d_{125c+37}(125n+37,62,87,112)&\equiv 0\pmod{5},\\
d_{125c+87}(125n+62,112)&\equiv 0\pmod{5},\\
d_{125c+92}(125n+82)&\equiv 0\pmod{5},
\end{aligned}
&
&
\begin{aligned}[t]
d_{125c+22}(125n+27,77,102)&\equiv 0\pmod{5},\\
d_{125c+47}(125n+52,77)&\equiv 0\pmod{5},\\
d_{125c+72}(125n+102)&\equiv 0\pmod{5},\\
d_{125c+104}(125n+95,120)&\equiv 0\pmod{5},\\
d_{125c+34}(125n+65,90,115)&\equiv 0\pmod{5},\\
d_{125c+59}(125n+40,90)&\equiv 0\pmod{5},\\
d_{125c+39}(125n+60,110)&\equiv 0\pmod{5},\\
d_{125c+114}(125n+35,60,85,110)&\equiv 0\pmod{5},\\
d_{125c+19}(125n+105)&\equiv 0\pmod{5},\\
d_{125c+24}(125n+100)&\equiv 0\pmod{5},\\
d_{125c+99}(125n+25)&\equiv 0\pmod{5},\\
d_{125c+124}(125n+50,75)&\equiv 0\pmod{5}.
\end{aligned}
\end{align*}
\end{theorem}
\begin{theorem}\label{thm15}
For all $c\geq 0$ and $n\geq 0$, 
\begin{align*}
&
\begin{aligned}[t]
d_{125c+76}(25n+23)&\equiv 0\pmod{25},\\
d_{125c+1}(125n+23,123)&\equiv 0\pmod{25},\\
d_{125c+106}(125n+93)&\equiv 0\pmod{25},\\
d_{125c+66}(125n+33,108)&\equiv 0\pmod{25},\\
d_{125c+91}(125n+108)&\equiv 0\pmod{25},\\
d_{125c+67}(125n+107)&\equiv 0\pmod{25},\\
d_{125c+103}(125n+96,121)&\equiv 0\pmod{25},\\
d_{125c+13}(125n+36,61,86,111)&\equiv 0\pmod{25},\\
d_{125c+63}(125n+61,111)&\equiv 0\pmod{25},\\
d_{125c+43}(125n+106)&\equiv 0\pmod{25},
\end{aligned}	
&
&
\begin{aligned}[t]
d_{125c+23}(125n+51,76)&\equiv 0\pmod{25},\\
d_{125c+48}(125n+101)&\equiv 0\pmod{25},\\
d_{125c+123}(125n+26,76,101)&\equiv 0\pmod{25},\\
d_{125c+99}(125n+75,100)&\equiv 0\pmod{25},\\
d_{125c+15}(125n+84)&\equiv 0\pmod{25},\\
d_{125c+115}(125n+109)&\equiv 0\pmod{25},\\
d_{125c+70}(125n+29,79,104)&\equiv 0\pmod{25},\\
d_{125c+95}(125n+54,79)&\equiv 0\pmod{25},\\
d_{125c+120}(125n+104)&\equiv 0\pmod{25}.
\end{aligned}
\end{align*}
\end{theorem}
\begin{theorem}\label{thm16}
For all $c\geq 0$ and $n\geq 0$, 
\begin{align*}
d_{125c}(125n+74,99,124)&\equiv 0\pmod{125}.
\end{align*}
\end{theorem}

\begin{remark}
The particular congruences $d_2(125n+97,122)\equiv 0\pmod{5}$ and $d_1(125n+23,123)\equiv 0\pmod{25}$ in Theorems \ref{thm14} and \ref{thm15} were first proved by Baruah, Das, and Talukdar \cite[Theorem 6.1]{bdt} using Radu's algorithm \cite{radu1}. Banerjee and Smoot \cite[Theorem 1.7]{bansmo} deduced the former congruences by showing their generalization given by
\begin{align*}
d_{125c+2}\left(5^{2\alpha+1}+5^{2\alpha}j+\dfrac{23\cdot 5^{2\alpha}+1}{8}\right)\equiv 0\pmod{5}
\end{align*}
for all $j\in\{1,2\}$, $\alpha\geq 1$, $c\geq 0$, and $n\geq 0$ using $q$-series techniques.

We organize the rest of the paper as follows. We enumerate in Section \ref{sec2} important identities necessary to prove our main results. These include the $5$-dissections of $f_1$ and $1/f_1$ involving the function 
\begin{align*}
R(q) := \prod_{n=1}^\infty\dfrac{(1-q^{5n-1})(1-q^{5n-4})}{(1-q^{5n-2})(1-q^{5n-3})}
\end{align*}
and certain formulas involving $R(q)$ and the parameter
\begin{align*}
K := \dfrac{f_2f_5^5}{qf_1f_{10}^5}
\end{align*}
due to Chern and Hirschhorn \cite{chehir}, and Chern and Tang \cite{chetn}. We present in Sections \ref{sec3}--\ref{sec6} the proofs of Theorems \ref{thm12} and \ref{thm14}--\ref{thm16} by utilizing these identities. We finally give in Section \ref{sec7} remarks about these proofs and pose conjectural congruences for $d_k(n)$ that need further investigation. We have done most of our computations using \textit{Mathematica}.
\end{remark}

\section{Important identities and $5$-dissections}\label{sec2}

We list in this section some helpful identities needed for the proofs of Theorems \ref{thm12} and \ref{thm14}--\ref{thm16}. We start with the following $5$-dissections \cite[(8.1.4), (8.4.4)]{hirsc}
\begin{align}
f_1 &= f_{25}\left(\dfrac{1}{R_5}-q-q^2R_5\right),\label{eq21}\\
\dfrac{1}{f_1} &= \dfrac{f_{25}^5}{f_5^6}\left(\dfrac{1}{R_5^4}+\dfrac{q}{R_5^3}+\dfrac{2q^2}{R_5^2}+\dfrac{3q^3}{R_5}+5q^4-3q^5R_5+2q^6R_5^2-q^7R_5^3+q^8R_5^4\right)\label{eq22},
\end{align}
where $R_m := R(q^m)$, and the following identities
\begin{align}
K+1 &= \dfrac{f_2^4f_5^2}{qf_1^2f_{10}^4},\label{eq23}\\
K-4 &= \dfrac{f_1^3f_5}{qf_2f_{10}^3},\label{eq24}
\end{align}
which are equivalent to \cite[(9.10)]{chehir} and \cite[(9.11)]{chehir}, respectively. For any $m\in\mathbb{N}$ and $n\in\mathbb{Z}$, we next define
\begin{align*}
P(m,n) := \dfrac{1}{q^mR_1^{m+2n}R_2^{2m-n}} + (-1)^{m+n}q^mR_1^{m+2n}R_2^{2m-n}.
\end{align*}
Chern and Tang \cite[Theorem 1.1]{chetn} obtained the recurrence formulas
\begin{align}
P(m,n+1) &= 4K^{-1}P(m,n) + P(m,n-1),\label{eq25}\\
P(m+2,n) &= KP(m+1,n)+P(m,n),\label{eq26}
\end{align}
with the initial values
\begin{align}
P(0,0) &= 2,\label{eq27}\\
P(0,1) &= 4K^{-1},\label{eq28}\\
P(1,0) &= K,\label{eq29}\\
P(1,-1) &= 4K^{-1}-2+K.\label{eq210}
\end{align}
We finally add the following congruence 
\begin{align*}
f_m^{5^k} \equiv f_{5m}^{5^{k-1}}\pmod{5^k}
\end{align*}
for all $k\geq 1$ and $m\geq 1$, which follows from the binomial theorem and will be frequently used without further notice.

\section{Proof of Theorem \ref{thm12}}\label{sec3}
In view of Theorem \ref{thm11}, we first prove the congruences
\begin{align}
d_{16}(25n+8)&\equiv 0\pmod{5},\label{eq31}\\
d_8(5n+1)&\equiv 0\pmod{25}.\label{eq32}
\end{align}

\begin{proof}[Proof of (\ref{eq31})]
We apply (\ref{eq21}) on the generating function for $d_{16}(n)$ so that
\begin{align}
\sum_{n\geq 0} d_{16}(n)q^n &= \dfrac{f_2^{16}}{f_1^{49}} \equiv \dfrac{f_{10}^{3}f_1f_2}{f_5^{10}}\nonumber\\
&\equiv \dfrac{f_{10}^{3}f_{25}f_{50}}{f_5^{10}}\left(\dfrac{1}{R_5}-q-q^2R_5\right)\left(\dfrac{1}{R_{10}}-q^2-q^4R_{10}\right)\pmod{5}.\label{eq33}
\end{align}
We look at the terms of (\ref{eq33}) involving $q^{5n+3}$, divide both sides by $q^3$, and then replace $q^5$ with $q$. We then have
\begin{align*}
\sum_{n\geq 0} d_{16}(5n+3)q^n &\equiv \dfrac{f_2^{3}f_5f_{10}}{f_1^{10}} \equiv \dfrac{f_{10}f_2^3}{f_5}\\
&\equiv\dfrac{f_{10}f_{50}^3}{f_5}\left(\dfrac{1}{R_{10}}-q^2-q^4R_{10}\right)^3 \pmod{5},
\end{align*}
where we apply (\ref{eq21}) on the last congruence. Considering the terms involving $q^{5n+1}$, dividing both sides by $q$, and then replace $q^5$ with $q$, we arrive at
\begin{align*}
\sum_{n\geq 0} d_{16}(25n+8)q^n &\equiv 5q\dfrac{f_2f_{10}^3}{f_1}\equiv 0\pmod{5}.
\end{align*}
\end{proof}

\begin{proof}[Proof of (\ref{eq32})]
Using (\ref{eq21}) on the generating function for $d_8(n)$ yields
\begin{align}
\sum_{n\geq 0} d_8(n)q^n &= \dfrac{f_2^{8}}{f_1^{25}} \equiv \dfrac{f_2^8}{f_5^5}\equiv \dfrac{f_{50}^8}{f_5^5}\left(\dfrac{1}{R_{10}}-q^2-q^4R_{10}\right)^8\pmod{25}\label{eq34}.
\end{align}	
We look at the terms of (\ref{eq34}) involving $q^{5n+1}$, divide both sides by $q$, and then replace $q^5$ with $q$, so that
\begin{align*}
\sum_{n\geq 0} d_8(5n+1)q^n &\equiv -125q^3\dfrac{f_{10}^8}{f_1^5}\equiv 0\pmod{25}.
\end{align*}
\end{proof}

To prove the last congruence of Theorem \ref{thm12}, it suffices to show that
\begin{align}
d_{25c+5}(25n+19)\equiv 0\pmod{25}\label{eq35},
\end{align}
as the case $k=0$ follows from the special case $p(25n+24)\equiv 0\pmod{25}$ of (\ref{eq12}) and the proof for the case $k=10$ is similar. We note that the case $c=0$ of (\ref{eq35}) follows from (\ref{eq11}). We also note that (\ref{eq35}) improves the congruence produced by combining the case $c=0$ and Theorem \ref{thm11}.

\begin{proof}[Proof of (\ref{eq35})]
We apply (\ref{eq21}) on the generating function for $d_{25c+5}(n)$ and obtain
\begin{align}
\sum_{n\geq 0} d_{25c+5}(n)q^n &= \dfrac{f_2^{25c+5}}{f_1^{75c+16}} \equiv \dfrac{f_{10}^{5c}f_1^9f_2^5}{f_5^{15c+5}}\nonumber\\
&\equiv \dfrac{f_{10}^{5c}f_{25}^9f_{50}^5}{f_5^{15c+5}}\left(\dfrac{1}{R_5}-q-q^2R_5\right)^9\left(\dfrac{1}{R_{10}}-q^2-q^4R_{10}\right)^5\pmod{25}.\label{eq36}
\end{align}
Extracting the terms of (\ref{eq36}) involving $q^{5n+4}$, dividing both sides by $q^4$, and then replacing $q^5$ with $q$, we get	
\begin{align}
\sum_{n\geq 0} d_{25c+5}(5n+4)q^n \equiv 5q^3\dfrac{f_2^{5c}f_5^9f_{10}^5}{f_1^{15c+5}}A\equiv5q^3\dfrac{f_{10}^{c+5}}{f_5^{3c-8}}A\pmod{25},\label{eq37}
\end{align}
where 
\begin{align*}
A := &-18P(3,1)-27P(3,2)+P(3,3)+73P(2,-1)+288P(2,0)+126P(2,1)\\
&-24P(2,2)+27P(2,3)+198P(1,2)+18P(1,-3)+126P(1,-2)+117P(1,-1)\\
&-234P(1,0)-378P(1,1)+81P(1,3)+2P(1,4)+9P(0,-4)-12P(0,-3)\\
&+252P(0,-2)+864P(0,-1)-803.
\end{align*}
Using the formulas (\ref{eq25})--(\ref{eq26}) and the initial values (\ref{eq27})--(\ref{eq210}), we find that
\begin{align}
A &= \dfrac{2816}{K^4}+\dfrac{9152}{K^3}+\dfrac{8976}{K^2}+\dfrac{660}{K}+485-264K+352K^2-44K^3\nonumber\\
&\equiv \dfrac{(1+K)^2(1+K^5)}{K^4}\equiv \dfrac{(1+K)^7}{K^4}\nonumber\\
&\equiv \left(\dfrac{f_2^4f_5^2}{qf_1^2f_{10}^4}\right)^7\left(\dfrac{qf_1f_{10}^5}{f_2f_5^5}\right)^4 \equiv \dfrac{f_2^{24}}{q^3f_1^{10}f_5^6f_{10}^8}\equiv \dfrac{f_2^4}{q^3f_5^8f_{10}^4}\pmod{5}.\label{eq38}
\end{align}
We combine (\ref{eq37}) and (\ref{eq38}) and use (\ref{eq21}) so that
\begin{align*}
\sum_{n\geq 0} d_{25c+5}(5n+4)q^n &\equiv 5q^3\dfrac{f_{10}^{c+5}}{f_5^{3c-8}}\cdot \dfrac{f_2^4}{q^3f_5^8f_{10}^4}\equiv 5\dfrac{f_{10}^{c+1}f_2^4}{f_5^{3c}}\\
&\equiv 5\dfrac{f_{10}^{c+1}f_{50}^4}{f_5^{3c}}\left(\dfrac{1}{R_{10}}-q^2-q^4R_{10}\right)^4 \pmod{25}.
\end{align*}
Considering the terms of the above congruence involving $q^{5n+3}$, dividing both sides by $q^3$, and then replacing $q^5$ with $q$, we finally arrive at
\begin{align*}
\sum_{n\geq 0} d_{25c+5}(25n+19)q^n &\equiv -25q\dfrac{f_2^{c+1}f_{10}^4}{f_1^{3c}}\equiv 0 \pmod{25}.
\end{align*}
\end{proof}

\section{Proof of Theorem \ref{thm14}}\label{sec4}
In this section, we prove Theorem \ref{thm14} by showing only the congruences
\begin{align}
d_6(125n+43,93)&\equiv 0\pmod{5},\label{eq41}\\
d_{61}(125n+38,63,88,113)&\equiv 0\pmod{5},\label{eq42}\\
d_{82}(125n+67,92,117)&\equiv 0\pmod{5},\label{eq43}
\end{align}
as the proofs for the remaining congruences are similar and follows from Theorem \ref{thm11}.

\begin{proof}[Proof of (\ref{eq41})]
We substitute (\ref{eq21}) to the generating function for $d_6(n)$ so that
\begin{align}
\sum_{n\geq 0} d_6(n)q^n &= \dfrac{f_2^6}{f_1^{19}} \equiv \dfrac{f_{10}f_1f_2}{f_5^4}\nonumber\\
&\equiv \dfrac{f_{10}f_{25}f_{50}}{f_5^4}\left(\dfrac{1}{R_5}-q-q^2R_5\right)\left(\dfrac{1}{R_{10}}-q^2-q^4R_{10}\right)\pmod{5}.\label{eq44}
\end{align}
We consider the terms of (\ref{eq44}) involving $q^{5n+3}$, divide both sides by $q^3$, and then replace $q^5$ with $q$. Applying (\ref{eq21}), we then get
\begin{align}
\sum_{n\geq 0} d_6(5n+3)q^n &\equiv \dfrac{f_2f_5f_{10}}{f_1^4} \equiv f_{10}f_1f_2\nonumber\\
&\equiv f_{10}f_{25}f_{50}\left(\dfrac{1}{R_5}-q-q^2R_5\right)\left(\dfrac{1}{R_{10}}-q^2-q^4R_{10}\right)\pmod{5}.\label{eq45}
\end{align}
We again consider the terms of (\ref{eq45}) involving $q^{5n+3}$, divide both sides by $q^3$, and then replace $q^5$ with $q$. We obtain 
\begin{align*}
\sum_{n\geq 0} d_6(25n+18)q^n &\equiv f_5f_{10}f_2 \equiv f_5f_{10}f_{50}\left(\dfrac{1}{R_{10}}-q^2-q^4R_{10}\right),\pmod{5}
\end{align*}
where the last congruence follows from (\ref{eq21}). By looking at the terms of the above congruence involving $q^{5n+1}$ and $q^{5n+3}$, we arrive at (\ref{eq41}).
\end{proof}

\begin{proof}[Proof of (\ref{eq42})]
We plug in (\ref{eq21}) to the generating function for $d_{61}(n)$ so that
\begin{align}
\sum_{n\geq 0} d_{61}(n)q^n &= \dfrac{f_2^{61}}{f_1^{184}} \equiv \dfrac{f_{10}^{12}f_1f_2}{f_5^{37}}\nonumber\\
&\equiv \dfrac{f_{10}^{12}f_{25}f_{50}}{f_5^{37}}\left(\dfrac{1}{R_5}-q-q^2R_5\right)\left(\dfrac{1}{R_{10}}-q^2-q^4R_{10}\right)\pmod{5}.\label{eq46}
\end{align}
We consider the terms of (\ref{eq46}) involving $q^{5n+3}$, divide both sides by $q^3$, and then replace $q^5$ with $q$. Utilizing (\ref{eq21}), we have
\begin{align}
\sum_{n\geq 0} d_{61}(5n+3)q^n &\equiv \dfrac{f_2^{12}f_5f_{10}}{f_1^{37}} \equiv \dfrac{f_{10}^3f_1^3f_2^2}{f_5^7}\nonumber\\
&\equiv \dfrac{f_{10}^3f_{25}^3f_{50}^2}{f_5^7}\left(\dfrac{1}{R_5}-q-q^2R_5\right)^3\left(\dfrac{1}{R_{10}}-q^2-q^4R_{10}\right)^2\pmod{5}.\label{eq47}
\end{align}
We extract the terms of (\ref{eq47}) involving $q^{5n+2}$, divide both sides by $q^2$, and then replace $q^5$ with $q$. We then get\\
\begin{align}
	\sum_{n\geq 0} d_{61}(25n+13)q^n &\equiv q\dfrac{f_2^3f_5^3f_{10}^2}{f_1^7}B \equiv qf_1^3f_2^3f_5f_{10}^2B\pmod{5},\label{eq48}
\end{align}
where, invoking (\ref{eq23}) and (\ref{eq25})--(\ref{eq210}),
\begin{align}
	B &:= -2P(1,1)-6P(0,1)-5= -\dfrac{32}{K}-9-2K\nonumber\\
	&\equiv \dfrac{3(K+1)^2}{K} \equiv 3\left(\dfrac{f_2^4f_5^2}{qf_1^2f_{10}^4}\right)^2\left(\dfrac{qf_1f_{10}^5}{f_2f_5^5}\right)\equiv \dfrac{3f_2^7}{qf_1^3f_5f_{10}^3}\pmod{5}.\label{eq49}
\end{align}
Combining (\ref{eq48}) and (\ref{eq49}) leads to 
\begin{align}
	\sum_{n\geq 0} d_{61}(25n+13)q^n &\equiv qf_1^3f_2^3f_5f_{10}^2\cdot \dfrac{3f_2^7}{qf_1^3f_5f_{10}^3} \equiv 3f_{10}\pmod{5}.\label{eq410}
\end{align}
Observe that the $q$-expansion of the right-hand side of (\ref{eq410}) contains only terms of the form $q^{10n}$. Hence, by collecting all the terms of (\ref{eq410}) involving $q^{5n+j}$ for $1\leq j\leq 4$, we get (\ref{eq42}).
\end{proof}

\begin{proof}[Proof of (\ref{eq43})]
We use (\ref{eq21}) on the generating function for $d_{82}$ so that
\begin{align}
	\sum_{n\geq 0} d_{82}(n)q^n &= \dfrac{f_2^{82}}{f_1^{247}} \equiv \dfrac{f_{10}^{16}f_1^3f_2^2}{f_5^{50}}\nonumber\\
	&\equiv \dfrac{f_{10}^{16}f_{25}^3f_{50}^2}{f_5^{50}}\left(\dfrac{1}{R_5}-q-q^2R_5\right)^3\left(\dfrac{1}{R_{10}}-q^2-q^4R_{10}\right)^2\pmod{5}.\label{eq411}
\end{align}
We extract the terms of (\ref{eq411}) involving $q^{5n+2}$, divide both sides by $q^2$, and then replace $q^5$ with $q$. Using (\ref{eq21}) yields
\begin{align}
	\sum_{n\geq 0} d_{82}(5n+2)q^n &\equiv q\dfrac{f_2^{16}f_5^3f_{10}^2}{f_1^{50}}B\equiv q\dfrac{f_{10}^5f_2}{f_5^7}\cdot\dfrac{3f_2^7}{qf_1^3f_5f_{10}^3}\equiv \dfrac{3f_{10}^3f_1^2f_2^3}{f_5^9}\nonumber\\
	&\equiv\dfrac{3f_{10}^3f_{25}^2f_{50}^3}{f_5^9}\left(\dfrac{1}{R_5}-q-q^2R_5\right)^2\left(\dfrac{1}{R_{10}}-q^2-q^4R_{10}\right)^3\pmod{5},\label{eq412}
\end{align}
where $B$ is given by (\ref{eq49}). By looking at the terms of (\ref{eq412}) involving $q^{5n+3}$, dividing both sides by $q^3$, and then replace $q^5$ with $q$, we deduce that
\begin{align}
	\sum_{n\geq 0} d_{82}(25n+17)q^n &\equiv 3q\dfrac{f_2^3f_5^2f_{10}^3}{f_1^9}C\equiv 3qf_{10}^3f_1f_2^3C\pmod{5},\label{eq413}
\end{align}
where, applying (\ref{eq23}) and (\ref{eq25})--(\ref{eq210}),
\begin{align}
	C &:= 2P(1,-1)+6P(1,0)-5= \dfrac{8}{K}-9+8K\nonumber\\
	&\equiv \dfrac{3(K+1)^2}{K} \equiv 3\left(\dfrac{f_2^4f_5^2}{qf_1^2f_{10}^4}\right)^2\left(\dfrac{qf_1f_{10}^5}{f_2f_5^5}\right)\equiv \dfrac{3f_2^7}{qf_1^3f_5f_{10}^3}\pmod{5}.\label{eq414}
\end{align}
Combining (\ref{eq413}) and (\ref{eq414}) yields 
\begin{align*}
	\sum_{n\geq 0} d_{82}(25n+17)q^n &\equiv 3qf_{10}^3f_1f_2^3\cdot \dfrac{3f_2^7}{qf_1^3f_5f_{10}^3} \equiv 9\dfrac{f_{10}^2f_1^3}{f_5^2}\pmod{5}\\
	&\equiv -\dfrac{f_{10}^2f_{25}^3}{f_5^2}\left(\dfrac{1}{R_5}-q-q^2R_5\right)^3\\
	&\equiv -\dfrac{f_{10}^2f_{25}^3}{f_5^2}\left(\dfrac{1}{R_5^3}-\dfrac{3q}{R_5^2}+5q^3-3q^5R_5^2-q^6R_5^3\right)\pmod{5},
\end{align*}
where we apply (\ref{eq21}) on the penultimate step. Taking all the terms of the above congruence involving $q^{5n+j}$ for $j\in \{2,3,4\}$, we get (\ref{eq43}).
\end{proof}

\section{Proof of Theorem \ref{thm15}}\label{sec5}
We establish Theorem \ref{thm15} by deriving only the congruences
\begin{align}
	d_{125c+76}(25n+23)&\equiv 0\pmod{25},\label{eq51}\\
	d_{125c+1}(125n+23,123)&\equiv 0\pmod{25},\label{eq52}
\end{align}
as the proofs of the other congruences are analogous.

\begin{proof}[Proof of (\ref{eq51})]
We apply (\ref{eq21}) and (\ref{eq22}) on the generating function for $d_{76}$ and get, modulo $25$,
\begin{align} 
\sum_{n\geq 0} d_{76}(n) &= \dfrac{f_2^{76}}{f_1^{229}}\equiv \dfrac{f_{10}^{15}f_2}{f_5^{45}f_1^4}\equiv \dfrac{f_{10}^{15}f_{50}}{f_5^{45}}\left(\dfrac{1}{R_{10}}-q^2-q^4R_{10}\right)\nonumber\\
&\times \dfrac{f_{25}^{20}}{f_5^{24}}\left(\dfrac{1}{R_5^4}+\dfrac{q}{R_5^3}+\dfrac{2q^2}{R_5^2}+\dfrac{3q^3}{R_5}+5q^4-3q^5R_5+2q^6R_5^2-q^7R_5^3+q^8R_5^4\right)^4.\label{eq53}
\end{align}
We extract the terms of (\ref{eq53}) involving $q^{5n+3}$, divide both sides by $q^3$, and then replace $q^5$ with $q$. We deduce that
\begin{align}
\sum_{n\geq 0} d_{76}(5n+3) &\equiv q^3\dfrac{f_2^{15}f_5^{20}f_{10}}{f_1^{69}}D \equiv q^3f_1^6f_2^{15}f_5^5f_{10}D\pmod{25},\label{eq54}
\end{align}
where
\begin{align}
D := &-4P(3,6)+40P(3,5)-418P(2,4)+1100P(2,3)-105P(2,5)\nonumber\\
&-1840P(1,2)+1200P(1,1)-1400P(1,3)-1500P(0,1)-1015.\label{eq55}
\end{align}
Set 
\begin{align*}
L :=\dfrac{f_1^3f_5}{qf_2f_{10}^3}
\end{align*}
so that $K=L+4$ by (\ref{eq24}). Using (\ref{eq24})--(\ref{eq210}) on (\ref{eq55}) yields
\begin{align}
D &:= -\dfrac{16384}{K^6}-\dfrac{91136}{K^5}-\dfrac{177664}{K^4}-\dfrac{172096}{K^3}-\dfrac{100864}{K^2}-\dfrac{38484}{K}-5711\nonumber\\
&+1986K+929K^2+36K^3\nonumber\\
&= K^{-6}(8\cdot 10^7L+112\cdot 10^6 L^2+682\cdot 10^5L^3+235\cdot 10^5L^4+4982500L^5\nonumber\\
&+659625L^6+52450L^7+2225L^8+36L^9)\nonumber\\
&\equiv 11K^{-6}L^9\equiv 11\left(\dfrac{f_1^3f_5}{qf_2f_{10}^3}\right)^9\left(\dfrac{qf_1f_{10}^5}{f_2f_5^5}\right)^6\equiv 11\dfrac{f_1^{33}f_{10}^3}{q^3f_2^{15}f_5^{21}}\pmod{25}.\label{eq56}
\end{align}
We deduce from (\ref{eq21}), (\ref{eq54}), and (\ref{eq56}) that
\begin{align}
	\sum_{n\geq 0} d_{76}(5n+3) &\equiv q^3f_1^6f_2^{15}f_5^5f_{10}\cdot 11\dfrac{f_1^{33}f_{10}^3}{q^3f_2^{15}f_5^{21}}\equiv 11\dfrac{f_1^{14}f_{10}^4}{f_5^{11}}\nonumber\\
	&\equiv 11\dfrac{f_{10}^4f_{25}^{14}}{f_5^{11}}\left(\dfrac{1}{R_5}-q-q^2R_5\right)^{14}\pmod{25}.\label{eq57}
\end{align}
We look for the terms of (\ref{eq57}) involving $q^{5n+4}$, divide both sides by $q^4$, and then replace $q^5$ with $q$, arriving at
\begin{align*}
	\sum_{n\geq 0} d_{76}(25n+23) &\equiv -11\cdot 5^6q^2\dfrac{f_2^4f_5^{14}}{f_1^{11}}\equiv 0\pmod{25}.
\end{align*}
Hence, we get $d_{76}(25n+23)\equiv 0\pmod{25}$ and (\ref{eq51}) follows from Theorem \ref{thm11}.
\end{proof}

\begin{proof}[Proof of (\ref{eq52})]
We begin with utilizing (\ref{eq21}) and (\ref{eq22}) on the generating function for $d_{125c+1}$ so that, modulo $25$,
\begin{align} 
	\sum_{n\geq 0} d_{125c+1}(n) &= \dfrac{f_2^{125c+1}}{f_1^{375c+4}}\equiv \dfrac{f_{10}^{25c}f_2}{f_5^{75c}f_1^4}\equiv \dfrac{f_{10}^{25c}f_{50}}{f_5^{75c}}\left(\dfrac{1}{R_{10}}-q^2-q^4R_{10}\right)\nonumber\\
	&\times \dfrac{f_{25}^{20}}{f_5^{24}}\left(\dfrac{1}{R_5^4}+\dfrac{q}{R_5^3}+\dfrac{2q^2}{R_5^2}+\dfrac{3q^3}{R_5}+5q^4-3q^5R_5+2q^6R_5^2-q^7R_5^3+q^8R_5^4\right)^4.\label{eq58}
\end{align}
We consider the terms of (\ref{eq58}) involving $q^{5n+3}$, divide both sides by $q^3$, and then replace $q^5$ with $q$. Using (\ref{eq21}), we then have, modulo $25$,
\begin{align}
	\sum_{n\geq 0} d_{125c+1}(5n+3)q^n &\equiv q^3\dfrac{f_2^{25c}f_5^{20}f_{10}}{f_1^{75c+24}}D\equiv  q^3\dfrac{f_2^{25c}f_5^{20}f_{10}}{f_1^{75c+24}}\cdot 11\dfrac{f_1^{33}f_{10}^3}{q^3f_2^{15}f_5^{21}}\nonumber\\
	&\equiv 11\dfrac{f_{10}^{5c-1}f_1^9f_2^{10}}{f_5^{15c+1}}\nonumber\\
	&\equiv 11\dfrac{f_{10}^{5c-1}f_{25}^9f_{50}^{10}}{f_5^{15c+1}}\left(\dfrac{1}{R_5}-q-q^2R_5\right)^9\left(\dfrac{1}{R_{10}}-q^2-q^4R_{10}\right)^{10},\label{eq59}
\end{align}
where $D$ is given by (\ref{eq56}). We collect the terms of (\ref{eq59}) involving $q^{5n+4}$, divide both sides by $q^3$, and then replace $q^5$ with $q$. This gives
\begin{align}
	\sum_{n\geq 0} d_{125c+1}(25n+23)q^n &\equiv 5q^5\dfrac{f_2^{5c-1}f_5^9f_{10}^{10}}{f_1^{15c+1}}E\equiv 5q^5\dfrac{f_{10}^{c+10}f_5^{9-3c}}{f_1f_2}E\pmod{25},\label{eq510}
\end{align}
where 
\begin{align*}
	E := &-18P(5,0)-54P(5,1)+7P(5,2)+73P(4,-2)+576P(4,-1)+882P(4,0)\\
	&+72P(4,1)+189P(4,2)+18P(3,-4)+252P(3,-3)+819P(3,-2)+702P(3,-1)\\
	&-2646P(3,0)-4104P(3,1)+486P(3,2)-108P(3,3)+18P(2,-5)-84P(2,-4)\\
	&-756P(2,-3)+6048P(2,-2)+16644P(2,-1)-5184P(2,0)-13608P(2,1)\\
	&-108P(2,2)-1386P(2,3)+6P(1,-5)+567P(1,-4)+4104P(1,-3)\\
	&-2268P(1,-2)-12636P(1,-1)-1053P(1,0)+19404P(1,1)+2358P(1,2)\\
	&-4158P(1,3)+9P(1,4)-162P(0,-4)+1296P(0,-3)+1134P(0,-2)\\
	&-44352P(0,-1)-9563.
\end{align*}
Applying (\ref{eq23}) and (\ref{eq25})--(\ref{eq210}) leads to
\begin{align}
	E &:=-\dfrac{12288}{K^5}-\dfrac{195072}{K^4}-\dfrac{27648}{K^3}-\dfrac{329152}{K^2}+\dfrac{512}{K}-2403-90747K-1448K^2\nonumber\\
	&-5502K^3+1712K^4-65K^5\nonumber\\
	&\equiv \dfrac{2(1+K)^9}{K^5}\equiv 2\left(\dfrac{f_2^4f_5^2}{qf_1^2f_{10}^4}\right)^9\left(\dfrac{qf_1f_{10}^5}{f_2f_5^5}\right)^5\equiv \dfrac{2f_1^2f_2}{q^4f_5^{10}f_{10}^5}.\pmod{5}\label{eq511}
\end{align}
We infer from (\ref{eq21}), (\ref{eq510}), and (\ref{eq511}) that
\begin{align*}
	\sum_{n\geq 0} d_{125c+1}(25n+23)q^n&\equiv 5q^5\dfrac{f_{10}^{c+10}f_5^{9-3c}}{f_1f_2}\cdot \dfrac{2f_1^2f_2}{q^4f_5^{10}f_{10}^5}\equiv 10q\dfrac{f_{10}^{c+5}f_1}{f_5^{3c+1}}\nonumber\\
	&\equiv 10q\dfrac{f_{10}^{c+5}f_{25}}{f_5^{3c+1}}\left(\dfrac{1}{R_5}-q-q^2R_5\right)\pmod{25}.
\end{align*}
Looking at the terms of the resulting congruence involving $q^{5n}$ and $q^{5n+4}$, we obtain (\ref{eq52}). 
\end{proof}

\section{Proof of Theorem \ref{thm16}}\label{sec6}

\begin{proof}
We first apply (\ref{eq22}) on the generating function for $d_{125c}$, yielding, modulo $125$,
\begin{align} 
	\sum_{n\geq 0} d_{125c}(n) &= \dfrac{f_2^{125c}}{f_1^{375c+1}}\equiv \dfrac{f_{10}^{25c}}{f_5^{75c}f_1}\nonumber\\
	&\equiv \dfrac{f_{10}^{25c}f_{25}^5}{f_5^{75c+6}}\left(\dfrac{1}{R_5^4}+\dfrac{q}{R_5^3}+\dfrac{2q^2}{R_5^2}+\dfrac{3q^3}{R_5}+5q^4-3q^5R_5+2q^6R_5^2-q^7R_5^3+q^8R_5^4\right).\label{eq61}
\end{align}
We look for the terms of (\ref{eq61}) involving $q^{5n+4}$, divide both sides by $q^4$, and then replace $q^5$ with $q$, so that, modulo $125$,
\begin{align} 
	\sum_{n\geq 0} d_{125c}(5n+4) &\equiv 5\dfrac{f_2^{25c}f_5^5}{f_1^{75c+6}}\equiv 5\dfrac{f_{10}^{5c}f_5^5}{f_5^{15c}f_1^6}\equiv 5\dfrac{f_{10}^{5c}f_{25}^{30}}{f_5^{15c+31}}\nonumber\\
	&\times\left(\dfrac{1}{R_5^4}+\dfrac{q}{R_5^3}+\dfrac{2q^2}{R_5^2}+\dfrac{3q^3}{R_5}+5q^4-3q^5R_5+2q^6R_5^2-q^7R_5^3+q^8R_5^4\right)^6, \label{eq62}
\end{align}
where we use (\ref{eq22}) on the last congruence. Collecting the terms (\ref{eq62}) involving $q^{5n+4}$, dividing both sides by $q^4$, and then replacing $q^5$ with $q$, we get.
\begin{align} 
	\sum_{n\geq 0} d_{125c}(25n+24) &\equiv 25q^4\dfrac{f_2^{5c}f_5^{30}}{f_1^{15c+31}}F\equiv 25q^4\dfrac{f_{10}^cf_5^{24-3c}}{f_1}F\pmod{125}.\label{eq63}
\end{align}
Using (\ref{eq24})--(\ref{eq210}), we find that 
\begin{align}
	F &:= 63P(4,0)+3728P(3,0)+27861P(2,0)+25404P(1,0)+106425\nonumber\\
	&= K^{-8}(63 K^{12}+4736 K^{11}+80913 K^{10}+574260 K^9+2441885 K^8+7183296 K^7\nonumber\\
	&+16217472 K^6+27033856 K^5+36552960 K^4+32194560 K^3+25591808 K^2\nonumber\\
	&+8257536 K+4128768)\nonumber\\
	&= K^{-8}(63 L^{12}+7760 L^{11}+355825 L^{10}+8865500 L^9+139368125 L^8+14887\cdot 10^5 L^7\nonumber\\
	&+112274\cdot 10^5 L^6+60764\cdot 10^6 L^5+23566\cdot 10^7 L^4+6416\cdot 10^8 L^3+1168\cdot 10^9 L^2\nonumber\\
	&+128\cdot 10^{10} L+64\cdot 10^{10})\nonumber\\
	&\equiv 3K^{-8}L^{12} \equiv 3\left(\dfrac{f_1^3f_5}{qf_2f_{10}^3}\right)^{12}\left(\dfrac{qf_1f_{10}^5}{f_2f_5^5}\right)^8 \equiv\dfrac{3f_1^{44}}{q^4f_5^{28}}\pmod{5}.\label{eq64}
\end{align}
We deduce from (\ref{eq21}), (\ref{eq63}), and (\ref{eq64}) that
\begin{align*} 
	\sum_{n\geq 0} d_{125c}(25n+24) &\equiv 25q^4\dfrac{f_{10}^cf_5^{24-3c}}{f_1}\cdot \dfrac{3f_1^{44}}{q^4f_5^{28}}\equiv 75\dfrac{f_{10}^cf_1^3}{f_5^{3c-4}}\\
	&\equiv 75\dfrac{f_{10}^cf_{25}^3}{f_5^{3c-4}}\left(\dfrac{1}{R_5}-q-q^2R_5\right)^3\\
	&\equiv 75\dfrac{f_{10}^cf_{25}^3}{f_5^{3c-4}}\left(\dfrac{1}{R_5^3}-\dfrac{3q}{R_5^2}+5q^3-3q^5R_5^2-q^6R_5^3\right)\pmod{125}.
\end{align*}
Taking all the terms involving $q^{5n+j}$ for $j\in \{2,3,4\}$ yields the desired congruence.
\end{proof}

\section{Closing remarks}\label{sec7}
We have shown in this paper infinite families of congruences modulo small powers of $5$ that are satisfied by $d_k(n)$, as provided by Theorems \ref{thm12} and \ref{thm14}--\ref{thm16}.
With that said, we remark that the list of these congruences is not yet exhaustive. In fact, numerical calculations via \textit{Mathematica} reveal the following conjectural congruences for $d_k(n)$ modulo small powers of $5$.

\begin{conjecture}\label{conj71}
For all $c\geq 0$ and $n\geq 0$, 
\begin{align*}
	d_{125c+58}(25n+16)&\equiv 0\pmod{125},\\
	d_{125c+83}(125n+41,91)&\equiv 0\pmod{125},\\
	d_{125c+100}(125n+124)&\equiv 0\pmod{125},\\
	d_{125c+5}(125n+69,119)&\equiv 0\pmod{125},\\
	d_{125c+30}(125n+69)&\equiv 0\pmod{125},\\
	d_{125c+60}(125n+14,64,89,114)&\equiv 0\pmod{125},\\
	d_{125c+58}(125n+91)&\equiv 0\pmod{625},\\
	d_{125c+58}(125n+66,116)&\equiv 0\pmod{3125}.
\end{align*}
\end{conjecture}

One may prove the particular congruences of Conjecture \ref{conj71} at $c=0$ using Radu's algorithm \cite{radu1} and/or Smoot's \texttt{RaduRK} package \cite{smoot1} and use these congruences to construct infinite congruences modulo arbitrary powers of $5$ via the localization method \cite{bansmo}, which heavily relies on modular forms. Nonetheless, to fully establish Conjecture \ref{conj71}, elementary proofs via $q$-series manipulations and dissections would be preferable. 

\section*{Acknowledgment}
The author would like to thank Pranjal Talukdar for bringing the recent paper of Yao \cite{yao} to his attention.

\end{document}